\documentclass[12pt]{article}
\usepackage{oldlfont}
\usepackage{enumerate}
\usepackage{amsmath}
\usepackage{amssymb}
\usepackage{amsthm}
\usepackage{mathrsfs}
\usepackage{amscd}
\usepackage{exscale}
\usepackage{latexsym}
\usepackage[all]{xy}
\usepackage{hyperref}
\usepackage{fancyhdr}
\usepackage{layout}
\usepackage{color}
\usepackage{graphicx}
\usepackage{amsfonts}
\usepackage{eucal}

\DeclareMathAlphabet{\mathpzc}{OT1}{pzc}{m}{it}

\begin{document}

\baselineskip=17pt

\pagestyle{headings}

\numberwithin{equation}{section}

\makeatletter                                                           

\def\section{\@startsection {section}{1}{\z@}{-5.5ex plus -.5ex         
minus -.2ex}{1ex plus .2ex}{\large \bf}}                                 


\pagestyle{fancy}
\renewcommand{\sectionmark}[1]{\markboth{ #1}{ #1}}
\renewcommand{\subsectionmark}[1]{\markright{ #1}}
\fancyhf{} 
\fancyhead[LE,RO]{\slshape\thepage}
\fancyhead[LO]{\slshape\rightmark}
\fancyhead[RE]{\slshape\leftmark}

\addtolength{\headheight}{0.5pt} 
\renewcommand{\headrulewidth}{0pt} 

\newtheorem{thm}{Theorem}
\newtheorem{mainthm}[thm]{Main Theorem}

\newcommand{\ZZ}{{\mathbb Z}}
\newcommand{\GG}{{\mathbb G}}
\newcommand{\Z}{{\mathbb Z}}
\newcommand{\RR}{{\mathbb R}}
\newcommand{\NN}{{\mathbb N}}
\newcommand{\GF}{{\rm GF}}
\newcommand{\QQ}{{\mathbb Q}}
\newcommand{\CC}{{\mathbb C}}
\newcommand{\FF}{{\mathbb F}}

\newtheorem{lem}[thm]{Lemma}
\newtheorem{cor}[thm]{Corollary}
\newtheorem{pro}[thm]{Proposition}
\newtheorem{proprieta}[thm]{Property}
\newcommand{\pf}{\noindent \textbf{Proof.} \ }
\newcommand{\eop}{$_{\Box}$  \relax}
\newtheorem{num}{equation}{}

\theoremstyle{definition}
\newtheorem{rem}[thm]{Remark}
\newtheorem*{D}{Definition}

\newcommand{\nsplit}{\cdot}
\newcommand{\G}{{\mathfrak g}}
\newcommand{\GL}{{\rm GL}}
\newcommand{\SL}{{\rm SL}}
\newcommand{\SP}{{\rm Sp}}
\newcommand{\LL}{{\rm L}}
\newcommand{\Ker}{{\rm Ker}}
\newcommand{\la}{\langle}
\newcommand{\ra}{\rangle}
\newcommand{\PSp}{{\rm PSp}}
\newcommand{\U}{{\rm U}}
\newcommand{\GU}{{\rm GU}}
\newcommand{\Aut}{{\rm Aut}}
\newcommand{\Alt}{{\rm Alt}}
\newcommand{\Sym}{{\rm Sym}}

\newcommand{\isom}{{\cong}}
\newcommand{\z}{{\zeta}}
\newcommand{\Gal}{{\rm Gal}}

\newcommand{\F}{{\mathbb F}}
\renewcommand{\O}{{\mathcal{O}}}
\renewcommand{\P}{{\mathcal{P}}}
\newcommand{\Q}{{\mathbb Q}}
\newcommand{\R}{{\mathbb R}}
\newcommand{\N}{{\mathbb N}}
\newcommand{\A}{{\mathcal{A}}}
\newcommand{\E}{{\mathcal{E}}}
\newcommand\topbot[2]{{\genfrac{}{}{0pt}{}{{#1}}{{#2}}}}

\newcommand{\DIM}{{\smallskip\noindent{\bf Proof.}\quad}}
\newcommand{\CVD}{\begin{flushright}$\square$\end{flushright}

\vskip 0.2cm\goodbreak}

\vskip 0.5cm

\title{On the local-global divisibility of torsion points on elliptic curves and ${\rm GL}_2$-type varieties}
\author{Florence Gillibert \footnote{Instituto de Matem\'aticas, Pontificia Universidad Cat\'olica de Valpara\'iso, Blanco Viel 596, Cerro Bar\'on, Valpara\'iso, Chile}, Gabriele Ranieri \footnote{Instituto de Matem\'aticas, Pontificia Universidad Cat\'olica de Valpara\'iso, Blanco Viel 596, Cerro Bar\'on, Valpara\'iso, Chile}}
\date{  }
\maketitle

\vskip 1.5cm

\begin{abstract}
Let $p$ be a prime number and let $k$ be a number field. 
Let $\E$ be an elliptic curve defined over $k$.
We prove that if $p$ is odd, then the local-global divisibility by any power of $p$ holds for the torsion points of $\E$.
We also show with an example that the hypothesis over $p$ is necessary.  

We get a weak generalization of the result on elliptic curves to the larger family of ${\rm GL}_2$-type varieties over $k$. In the special case of the abelian surfaces $A/k$ with quaternionic multiplication over $k$ we obtain that for all prime $p$, except a finite number depending on $A$, the local-global divisibility by any power of $p$ holds for the torsion points of $A$. 

\end{abstract}

\section{Introduction}

Let $k$ be a number field and let ${\mathcal{A}}$ be a commutative algebraic group defined over $ k $.
Several papers have been written on  the following classical question, known as \emph{Local-Global Divisibility Problem}.

\par\bigskip\noindent  {\sc Problem}: \emph{Let $P \in {\mathcal{A}}( k )$ and $q$ be a positive integer. Assume that for all but finitely many valuations $v$ of $k$, there exists $D_v \in {\mathcal{A}}( k_v )$ such that $P = qD_v$. Is it possible to conclude that there exists $D\in {\mathcal{A}}( k )$ such that $P=qD$?}

\par\bigskip\noindent  By  B\'{e}zout's identity, to get answers for a general integer it is sufficient to solve it for powers $p^n$ of a prime. In the classical case of ${\mathcal{A}}={\mathbb{G}}_m$, the answer is positive for any number field for $ p $ odd, and negative for instance for $k=\Q$ and $q=8$ (and $P=16$) (see for example \cite{AT}, \cite{T}).

\bigskip  For general commutative algebraic groups, Dvornicich and Zannier gave a cohomological  interpretation of the problem (see \cite{DZ}, \cite{DZ2} and \cite{DZ3}) that we shall explain.
Let $ \Gamma $ be a group and let $ M $ be a $\Gamma$-module.
We say that a cocycle $Z \colon \Gamma \rightarrow M$ satisfies the local conditions if for every $\gamma \in \Gamma$ there exists $m_\gamma \in M$ such that $Z_\gamma = \gamma ( m_\gamma ) - m_\gamma$.
The set of the classes of the cocycles in $ H^1 ( \Gamma , M ) $ that satisfy the local conditions is a subgroup of $H^1 ( \Gamma, M )$.
We call it the first local cohomology group $H^1_{{\rm loc}} ( \Gamma, M )$.
Equivalently,
\[
H^1_{{\rm loc}} ( \Gamma , M ) =  \cap_{ C \leq \Gamma } \ker ( H^1 ( \Gamma , M ) \rightarrow H^1 ( C , M ) ),
\]
where $ C $ varies among the cyclic subgroups of $ \Gamma $ and the above maps are the restrictions.   
Dvornicich and Zannier \cite[Proposion 2.1]{DZ} proved the following result.

\begin{pro}\label{pro1}
Let $ p $ be a prime number, let $n$ be a positive integer, let $ k $ be a number field and let $ \A $ be a commutative algebraic group defined over $ k $.
If $H^1_{{\rm loc}} ( \Gal ( k ( \A[p^n] ) / k ) , \A[p^n] ) = 0$, then the local-global divisibility by $p^n$ over $\A ( k )$ holds.
\end{pro}

The converse of Proposition \ref{pro1} is not true. However, in the case when the group $H^1_{{\rm loc}} ( \Gal ( k ( \A[p^n] ) / k ) , \A[p^n] )$ is not trivial, we can find an extension $ L $ of $ k $ such that $L$ and $k ( \A[p^n] )$ are linearly disjoint over $ k $, and such that the local-global divisibility by $ p^n $ over $ \A ( L ) $ does not hold (see \cite[Theorem 3]{DZ3} for the details).

The local-global divisibility problem has been studied in two interesting families of commutative algebraic groups: the algebraic tori (see \cite{DZ} and \cite{I}) and the elliptic curves (see \cite{DZ}, \cite{DZ2}, \cite{DZ3}, \cite{P}, \cite{PRV}, \cite{PRV2}).
Let $p$ be a prime number, let $ k $ be a number field and let $ \E $ be an elliptic curve defined over $k$.  
Dvornicich and Zannier (see \cite[Theorem 1]{DZ3}) found a criterion for the validity of the local-global divisibility by a power of $ p $ over $ \E ( k ) $.
In \cite{PRV} and \cite{PRV2} Paladino, Viada and the second author refined this criterion.  
In particular they proved that for every positive integer $r$, there exists a constant $C ( r )$, only depending on $ r $, such that for every number field $ k $ of degree $\leq r$, for every elliptic curve $\E$ defined over $ k $, for every prime number $p > C ( r )$, and for every $n \in \N$, the local-global divisibility by $p^n$ holds for $\E ( k )$ (see \cite[Theorem 1 , Corollary 2]{PRV}).  As we shows in this article (cf. our Theorem \ref{teo2} below) if we restrict our study to the torsion points of an elliptic curves $ \E $ defined over a number field $k$ we obtain a stronger result. In fact we show that for any odd prime $p$ and any $n\in \N$, the local-global divisibility by $p^n$ holds for the torsion points of $ \E ( k ) $. 

In this paper, we focus on a particular family of abelian varieties (the ${\rm GL}_2$-type varieties), but we restrict our study of the validity of the local-global divisibility principle to the torsion points.
Elliptic curves are an important particular case of ${\rm GL}_2$-type varietes, and so we devote a section to study it.
More precisely, we prove the following result. 

\begin{thm}\label{teo2}
Let $ \E $ be an elliptic curve defined over a number field $k$. Let $ p $ be an odd prime number. 
Then the local-global divisibility by any power of $p$ holds for the torsion points of $\E$.
\end{thm}

An abelian variety defined over a number field $ k $ is said to be of ${\rm GL}_2$-type if there is an embedding $\phi \colon E \hookrightarrow {\rm End}_k ( \A ) \otimes \Q$, where $E$ is a number field such that $[E: \Q] = {\rm dim} ( \A )$ (see \cite[Chapter 2, Chapter 5]{R}). 
Let $ \O_E $ be the ring of the algebraic integers of $E$.
Then $R = E \cap \phi^{-1} \left( {\rm End}_k ( \A ) \otimes \Z \right)$ is an order of $\O_E$ which is isomorphic to some subring of ${\rm End}_k ( \A )$ via $\phi$.
Following \cite[Chapter 2]{R}, we say that a prime number is {\bf good} for $ \A $ if it does not divide the index $[\O_E : R]$ (the notion of good depends also on $E$). 
We prove the following result, generalizing the result on elliptic curves.

\begin{thm}\label{teo3}   
Let $ \A $ be a ${\rm GL}_2$-type variety defined over a number field $k$. Let $E$ be a number field of degree ${\rm dim} ( \A )$ that embeds into ${\rm End}_k ( \A ) \otimes \Q$. Let $p$ be a {\bf good} odd prime number, such that $p$ splits completely in $E$. Then the local-global divisibility by any power of $p$ holds for the torsion points of $\A$.
\end{thm}  

Observe that for every number field $F$ of degree $\geq 2$ there are infinitely many prime numbers $p$ such that $( p )$ does not split completely in $F$.
Then, if only one number field $E$ of degree equal to to the dimension of $\A$ embeds into ${\rm End}_k ( \A ) \otimes \Q$, to apply Theorem \ref{teo3} we need to exclude infinitely many prime numbers.
We now study a case where the number field $E$ is not unique and we can exclude only a finite number of prime numbers.

Suppose that ${\rm End}_k ( \A )$ is an order in a quaternion algebra.
Then we say that $\A$ is an abelian variety with quaternionic multiplication over $k$.
In the special case where $\A$ is of dimension $2$ and it has quaternionic multiplication over $k$ by an order $\O$ in an indefinite quaternion algebra $B/ \Q$, we have the following result. %

\begin{cor}\label{cor4}
Let $B$ be an indefinite quaternion algebra over $\Q$. Let $\A/k$ be an abelian variety of dimension $2$ such that ${\rm End}_k ( \A )$ is isomorphic to an order $\O$ of $B$. Let $\O_m$ be a maximal order of $B$ such that $\O\subseteq \O_m$. Let $p\neq 2$ be a prime number such that $p$ does not divide the discriminant $D$ of the quaternion algebra $B$, and $p$ does not divide the index $M=[\O_m:\O]$. For any $n\geq 1$ the local-global divisibility by $p^n$ holds for torsion points of $\A$. 
\end{cor}
For abelian varieties $\A/k$ of higher dimension with quaternionic multiplication, we would exclude infinitely many primes when applying Theorem $\ref{teo3}$. In fact $A$ would have multiplication by a quaternion algebra $B$ over some number field $F$, and the condition $F\subseteq E$ would just remove the infinitely many primes numbers that do not split in $F$.

Here is the plan of the paper.

Let $ p $ be a prime number and let $n$ be a positive integer. 
Let $ \A $ be a commutative connected algebraic group defined over a number field $k$.
In Section \ref{sec1} we prove that if the local-global divisibility by $p^n$ does not hold for $\A$ over $ k $ and it fails for a torsion point, then $\A$ has a torsion point of order $p^n$ over $k$.

In Section \ref{sec2} we use this to make some calculations of the first local cohomology group that allow us to prove Theorem \ref{teo2} in section \ref{sec3}.
In this section we also show with an example that we need the hypothesis $p$ odd in Theorem \ref{teo2}.

In Section \ref{sec4} we prove Theorem \ref{teo3}.

Finally in Section \ref{sec5} we prove Corollary \ref{cor4}.   
\vspace{5 pt}

{\bf Aknowledgments.} The first author was supported by the project Fondecyt Iniciaci\'on Number 11130409.  The second author was supported by the project Fondecyt Regular Number 1140946.

\section{Preliminary results}\label{sec1}

Let $ p $ be a prime number, let $m$ be a positive integer and let $ \A $ be a connected commutative algebraic group defined over a number field $k$.

\begin{lem}\label{lem11}
Suppose that there exists $P \in \A ( k )$ a torsion point such that $ P $ is divisible by $p^m$ over $\A ( k_v )$ for all but finitely many valuations $ v $ of $ k $ and $P$ is not divisible by $p^m$ over $\A ( k )$. Then there exists $Q \in \A ( k )$ such that $Q$ has order equal to the power of $p$ dividing the order of $P$, and the local-global divisibility by $p^m$ fails for $Q$ over $\A ( k )$.
\end{lem}

\DIM Suppose that $P$ has order $l p^a$ with $l$ not divided by $p$ and $a \in \N$.
Then every multiple of $ P $, and so also $lP$, is divisible by $p^m$ for all but finitely many $ v $ over $\A ( k_v )$.
Let $D \in \A ( \overline{k} )$ be such that $P = p^m D$.
Suppose that $l P$ is divisible by $p^m$ in $\A ( k )$. 
Then there exists $R \in \A ( k )$ such that $l P = p^m R$.
But since $P = p^m D$, we have $l P = p^m l D$.
Thus $p^m ( R - l D ) = 0$, which gives
\[
l D = R + C
\]
with $C \in \A[p^m]$.
Since $ p $ does not divide $l$, there exists $C^\prime \in \A[p^m]$ such that $l C^\prime = C$.
Then $l ( D - C^\prime ) = R$.
Since $R \in \A ( k )$, we have $l ( D - C^\prime ) \in \A ( k )$.
Moreover $p^m ( D - C^\prime ) = P$ because $C^\prime \in \A[p^m]$.
Thus $p^m ( D - C^\prime ) \in \A ( k )$.
Since $( l, p^m ) = 1$, by B\'ezout identity, we then get $D - C^\prime \in \A ( k )$.
Moreover we have already remarked that $p^m ( D - C^\prime ) = P$.
Then $P$ is divisible by $p^m$ in $\A ( k )$, which is a contradiction
\CVD 

The following lemma gives a strong condition on the Galois action over the $p$-torsion part of $\A$.
Using this condition, we shall calculate the first local cohomology group of $ \A[p^m] $ as module over $\Gal ( k ( \A[p^m] ) / k )$.

\begin{lem}\label{lem12}
Let $P \in \A ( k )$ be a torsion point such that the local-global divisibility by a power of $p$ fails for $P$ over $\A ( k )$ and let $p^n$ be the minimal power of $p$ for which this occurs.
Then there exists $Q \in \A ( k )$ such that $p^{n-1} Q = P$. In particular, $\A$ admits a point of order $p^n$ defined over $k$. 
\end{lem}     

\DIM By Lemma \ref{lem11} and the fact that $0$ is divisible by any power of $p$, $P$ has order divided by $p^r$, with $r \geq 1$. 
Since $p^n$ is minimal, there exists $Q \in \A ( k )$ such that $p^{n-1} Q = P$.
In particular $Q$ has order divided by $p^{n + r - 1} \geq p^n$.
Then $\A ( k )$ has a $k$-rational point of order $p^n$.
\CVD
 
\section{Computing the first local cohomology group}\label{sec2}

Let $p$ be a prime number. For every positive integer $n$ set $V_2 = ( \Z / p^n \Z )^2$.

\begin{pro}\label{pro21}
Suppose that $p$ is odd. Let $n \in \N^\ast$ and let $G$ be a subgroup of ${\rm GL}_2 ( \Z / p^n \Z )$ such that there exists $Q \in V_2$ of order $p^n$ such that for all $g \in G$, $g ( Q ) = Q$. Then $H^1_{{\rm loc}} ( G , V_2 ) = 0$.
\end{pro}

\DIM Since $G$ fixes an element of $V_2$ of order $p^n$, we can fix a basis of $V_2$ such that every element of $G$ can be written as  
\[
\left(
\begin{array}{cc}
1 & f \\
0 & c \\
\end{array}
\right)    
\]
for certain $f \in \Z / p^n \Z$, $c \in ( \Z / p^n \Z )^\ast$.
Observe that $G$ has a unique $p$-Sylow $H$, given by the elements of determinant congruent to $1$ modulo $p$.
Since $V_2$ is a $p$-group, the restriction $H^1 ( G , V_2 ) \rightarrow H^1 ( H , V_2 )$ is injective and it induces a map of the first local cohomology groups. 
Then we can suppose that $G$ is a $p$-group. 
Thus for every element of $G$, there exist $f \in \Z / p^n \Z, s \in ( \Z / p^n \Z )^\ast$ and $k \in \N^\ast$, such that such an element can be written
\[     
\left(
\begin{array}{cc}
1 & f \\
0 & 1 + s p^k \\
\end{array}
\right).    
\]
Since the first local cohomology group of a module over a cyclic group is always trivial, we can suppose that $G$ has at least two generators.
In the following lemma we give a set of generators for $G$.

\begin{lem}\label{lem22}
Let $i$ be the biggest integer such that, for every
\[
\tau =
\left(
\begin{array}{cc}
1 & f \\
0 & 1 + s p^k \\
\end{array}
\right) \in G,    
\]
we have $k \geq i$, and let $e \in \N$ be such that
\[
\left(
\begin{array}{cc}
1 & p^e \\
0 & 1 \\
\end{array}
\right)    
\]
generates the group of the upper unitriangular matrices of $G$.
Then there exists $b \in \Z / p^n \Z$ such that $G$ is generated by the two elements
\[
\delta =
\left(
\begin{array}{cc}
1 & b \\
0 & 1 + p^i  \\
\end{array}
\right), \    
\sigma =
\left(
\begin{array}{cc}
1 & p^e \\
0 & 1 \\
\end{array}
\right).    
\]

Finally, if $p^e$ divides $b$, we can replace $\delta$ with 
\[
\left(
\begin{array}{cc}
1 & 0 \\
0 & 1 + p^i  \\
\end{array}
\right).   
\]
\end{lem}

\DIM By definition of $i$, there exist $f \in \Z / p^n \Z$ and $s \in ( \Z / p^n \Z )^\ast$ such that 
\[
\gamma =  
\left(
\begin{array}{cc}
1 & f \\
0 & 1 + s p^i \\
\end{array}
\right) \in G.     
\]
Since $p \geq 3$ and $s$ is invertible, $( 1 + s p^i )$ generates the subgroup of the elements of $( \Z / p^n \Z )^\ast$ congruent to $1$ modulo $p^i$. 
Then, there exists $g \in \N$ such that $( 1 + s p^i )^g \equiv 1 + p^i \mod ( p^n )$. 
Let $\delta = \gamma^g$.
We have
\[
\gamma^g =
\left(
\begin{array}{cc}
1 & b \\
0 & 1 + p^i 
\end{array}
\right),  
\]
for one $b \in \Z / p^n \Z$.
Let 
\[
\tau =
\left(
\begin{array}{cc}
1 & f^\prime \\
0 & 1 + s^\prime p^k \\
\end{array}
\right) \in G.
\]
Since $k \geq i$, there exists $j \in \N$ such that
\[
( 1 + p^i )^j \equiv s^\prime p^k + 1 \mod ( p^n ).
\]
Then $\tau \delta^{-j}$ is an upper unitriangular matrix.
By definition $\sigma$ is the generator of the group of the upper unitriangular matrices in $G$.
Thus there exists an integer $c$ such that $\tau \delta^{-j} = \sigma^c$, which gives that $\tau$ is in the group generated by $\delta$ and $\sigma$.
Now observe that if $\sigma$ is the identity, then the group has just a generator, that is contrary with the hypothesis.

Suppose that $p^e$ divides $b$.
Let $\widetilde{b} \in \Z$ be an integer in the class of $b$ and let $h \in \Z / p^n \Z$ be such that $h$ is the class of $- \widetilde{b} / p^e$.  
Then $\delta \sigma^h$ is a diagonal matrix.
Moreover $\{ \delta \sigma^h, \sigma \}$ still generate $G$  
because $\delta \in \langle \delta \sigma^h, \ \sigma \rangle$.     
\CVD

By Lemma \ref{lem22} we have two cases: either $\delta$ is diagonal or $\delta$ is not diagonal.
We first prove that $H^1_{{\rm loc}} ( G, V_2 ) = 0$ in the case where $\delta$ is diagonal.
Then we use this to study the other case.
\bigskip

{\bf The case where $\delta$ is diagonal}. Suppose that $G = \langle \delta, \sigma \rangle$ with
\[
\delta =
\left(
\begin{array}{cc}
1 & 0 \\
0 & 1 + p^i  \\
\end{array}
\right), \    
\sigma =
\left(
\begin{array}{cc}
1 & p^e \\
0 & 1 \\
\end{array}
\right).    
\]   
Let $Z$ be a cocycle representing a class in $H^1_{{\rm loc}} ( G, V_2 )$.
Then there exist $v, w \in V_2$ such that $Z_\delta = \delta ( v ) - v$, $Z_\sigma = \sigma ( w ) - w$.
By summing to $Z$ the coboundary $Z^\prime_\tau = \tau ( ( - v ) ) - ( - v )$ for every $\tau \in G$, we get a new cocycle in the class of $Z$, which is $0$ over $\delta$. 
Then we can suppose $Z_\delta = ( 0, 0 )$.
Moreover for every $( \alpha, \beta ) \in ( \Z/ p^n \Z )^2$, we have $( \sigma - Id ) ( ( \alpha, \beta ) ) = ( p^e \beta, 0 )$.
Thus we can suppose 
\[
Z_\delta = ( 0, 0 ), \ Z_\sigma = ( p^e \beta, 0 )
\]
for a certain $\beta \in \Z / p^n \Z$.
Since $Z$ is a cocycle, we have $Z_{\delta \sigma} = Z_\delta + \delta Z_\sigma$, which is equal to $( p^e \beta, 0 )$.
Since the class of $Z$ is in $H^1_{{\rm loc}} ( G, V_2 )$, we have $Z_{\delta \sigma} = ( \delta \sigma - Id ) ( t )$ for a certain $t \in V_2$.
But 
\[
\delta \sigma =
\left(
\begin{array}{cc}
1 & p^e \\
0 & 1 + p^i \\
\end{array}
\right).             
\]
Then the system $( \delta \sigma - Id ) ( ( x, y ) ) = ( p^e \beta, 0 )$ gives 
\[
\begin{cases}
p^e y\equiv p^e \beta  \mod ( p^n )\\
p^i y  \equiv 0  \mod ( p^n )
\end{cases}
\]
and it should have a solution.
Suppose first $e \geq i$. 
Then $p^e y \equiv 0 \mod ( p^n )$.
Thus $p^e \beta \equiv 0 \mod ( p^n )$.
Then $Z_\sigma = ( 0, 0 )$.
Since $Z_\delta = ( 0, 0 )$, $Z_\sigma = ( 0, 0 )$ and $\sigma$, $\delta$ generate $G$ by hypothesis, for every $\tau \in G$, $Z_\tau = ( 0, 0 )$.
Thus $Z$ is a coboundary.

Suppose now $e < i$.
Since $p^e y \equiv p^e \beta \mod ( p^n )$, we have $p^i y \equiv p^i \beta \mod ( p^n )$ and $p^i \beta \equiv 0 \mod ( p^n )$.
Thus $Z_\delta = ( \delta - Id ) ( ( 0, \beta ) )$ and $Z_\sigma = ( \sigma - Id ) ( ( 0, \beta ) )$.
But by Lemma \ref{lem22}, $\delta$ and $\sigma$ generate $G$.
Thus for every $\tau \in G$, $Z_\tau = ( \tau - Id ) ( ( 0, \beta ) )$.
Then $Z$ is a coboundary.
We have proved that if $\delta$ is diagonal, then $H^1_{{\rm loc}} ( G, V_2 ) = 0$.
\bigskip

{\bf The case when $\delta$ is not diagonal}. We consider now the case in which $\delta$ is not diagonal.
Consider the subgroup $\widetilde{G}$ of $G$ of the $\tau \in G$ such that
\[
\tau =
\left(
\begin{array}{cc}
1 & b^\prime \\
0 & 1 + s p^k  
\end{array}
\right),
\]
with $s \in \Z / p^n \Z^\ast$, $k \in \N$ and $p^e$ dividing $b^\prime$.
Then by Lemma \ref{lem22}, $\widetilde{G}$ is generated by a diagonal matrix $\delta^\prime$ and by $\sigma$.
The previous case shows that $H^1_{{\rm loc}} ( \widetilde{G} , V_2 ) = 0$.
Then, if $Z$ is a cocycle representing a class in $H^1_{{\rm loc}} ( G, V_2 )$, by summing a coboundary we can suppose that $Z_\tau = ( 0, 0 )$ for every $\tau \in \widetilde{G}$.
Write 
\[
\delta =
\left(
\begin{array}{cc}
1 & p^l m \\
0 & 1 + p^i  \\
\end{array}
\right)
\]
with $l < e$ (if not we are in the previous case by Lemma \ref{lem22}), $m \in ( \Z / p^n \Z )^\ast$. 
Consider $( \alpha , \beta ) \in V_2$ such that $Z_\delta = \delta ( ( \alpha , \beta ) ) - ( \alpha , \beta )$.
Since $p$ is odd, a short computation gives
\[
\delta^{p^{e-l}} =
\left(
\begin{array}{cc}
1 & p^e m^\prime \\
0 & ( 1 + p^i )^{p^{e-l}}  \\
\end{array}
\right)
\]
with $m^\prime \in ( \Z / p^n \Z )^\ast$.
Then $\delta^{p^{e-l}} \in \widetilde{G}$.
Thus 
\[
Z_{\delta^{p^{e-l}}} = \delta^{p^{e-l}} ( ( \alpha , \beta ) ) - ( \alpha , \beta ) = ( 0 , 0 ).
\]
Thus $p^e \beta \equiv 0 \mod ( p^n )$ and so $Z_\sigma = \sigma ( ( \alpha , \beta ) ) - ( \alpha , \beta )$.
Then, since $\delta$ and $\sigma$ generate $G$ by Lemma \ref{lem22} and $Z_\delta = \delta ( ( \alpha , \beta ) ) - ( \alpha , \beta )$, $Z_\sigma = \sigma ( ( \alpha , \beta ) ) - ( \alpha , \beta )$, we get that $Z$ is a coboundary.
\CVD     
 
\section{The case of elliptic curves}\label{sec3}

In this section we prove Theorem \ref{teo2}.
Then we show with a counterexample that the theorem is false in the case $p = 2$.
\bigskip

{\bf Proof of Theorem \ref{teo2}}. Suppose that there exists $P \in \E ( k )$ a torsion point for which the local-global divisibility by a power of $p$ fails over $k$.
By Lemma \ref{lem11} we can suppose that $P$ has order a power of $p$.
Let $n$ be the minimal integer such that $P$ is divisible by $p^n$ over $\E ( k_v )$ for all but finitely many primes $v$ and $P$ is not divisible by $p^n$ over $\E ( k )$.
Then, by Lemma \ref{lem12}, $\E$ admits a $k$-rational point of order $p^n$.
Conclude by applying Propositions \ref{pro21} and \ref{pro1}.
\CVD

To construct a counterexample to Theorem \ref{teo2} in the case $p=2$, we first recall the following result that we shall use also in the following section.

\begin{pro}\label{profin}
Let $\A$ be an abelian variety defined over a number field $k$, let $p$ be a prime number and let $n$ be a positive integer.  
Suppose that there exists a cocycle $Z$ representing a class $[Z]$ in $H^1_{{\rm loc}} ( \Gal ( k ( \A[p^n] ) / k ) , \A[p^n] )$, and $[Z] \neq 0$.
Moreover suppose that there exists $D \in \A ( K( \A[p^{n}] ) )$ such that for every $\tau \in \Gal ( k ( \A[p^n] ) / k )$, we have $\tau ( D ) = D + Z_{\tau}$. Then $P = p^n D$ is divisible by $p^n$ in $\A ( k_v )$ for every prime $v$ of $k$ not ramified over $k ( \A[p^n])$, but $P$ is not divisible by $p^n$ over $\A ( k )$.
\end{pro}

\DIM It is well-known. For a proof, see for instance the proof of \cite[Remark 2.2]{DZ}, or reverse the proof of \cite[Proposition 2.1]{DZ}.
\CVD                    
 
The following corollary, that could be interesting also in more general cases, will be applied to find a counterexample to Theorem \ref{teo2} in the case $p = 2$.

\begin{cor}\label{cor31}  
Let $n$ be a positive integer, let $r$ be an odd positive integer, let $p$ be a prime number and let $G$ be a subgroup of ${\rm GL}_r ( \Z / p^n \Z )$ acting on $( \Z / p^n \Z )^r$. Suppose that there exists $Z \colon G \rightarrow ( \Z / p^n \Z )^r$ a cocycle representing a class $[Z]$ in $H^1_{{\rm loc}} ( G , ( \Z / p^n \Z )^r )$, such that $[Z]\neq 0$. 
Suppose also that there exists an abelian variety $\A$ of dimension $d = ( r+1 )/ 2$ defined over a number field $k$, such that: 
\begin{enumerate}
\item There exist $P_1, P_2, \ldots , P_r \in \A[p^n]$ independent over $\Z / p^n \Z$, an isomorphism $\phi \colon \Gal ( k ( \A[p^n] / k ) \rightarrow G$ and a $\Z / p^n \Z$-linear isomorphism $\tau \colon \{ P_1 , P_2 , \ldots , P_r \} \rightarrow ( \Z / p^n \Z )^r$, such that for all $\sigma \in \Gal ( k ( \A[p^n] / k )$ and $v \in \{ P_1 , P_2 , \ldots , P_r \}$, 
\[
\phi ( \sigma ) \cdot \tau ( v ) = \tau ( \sigma ( v ) );
\]
\item There exists $Q \in \A ( k ( \A[p^{n}] ) )$ such that for all $\sigma \in {\rm Gal} ( k ( \A[p^n]) / k )$, we have $\sigma ( Q ) = Q + \tau^{-1} ( Z_{\phi ( \sigma )} )$.
\end{enumerate}
Then $p^n Q$ is divisible by $p^n$ over $\A ( k_v )$ for all prime $v$ not ramified over $k ( \A[p^n])$, but it is not divisible by $p^n$ over $\A ( k )$.
\end{cor}

\DIM Define $W \colon \Gal ( k ( \A[p^n] ) / k ) \rightarrow \A[p^n]$ as $W = \tau^{-1} \circ Z \circ \phi$.
Since $Z$ is a cocycle, a simple verification shows that $W$ is a cocycle. 
Let $\sigma$ be in ${\rm Gal} ( k ( \A[p^n] ) / k )$ and let $a$ be in $\Z / p^n \Z$ such that $Z_{\phi ( \sigma )} = \phi ( \sigma ) \cdot a - a$.
Then $W_\sigma = \sigma ( \tau^{-1} ( a ) ) - \tau^{-1} ( a )$, and so $W$ satisfies the local conditions.
Finally observe that if $W$ were a coboundary, then also $Z$ should be a coboundary and so $W$ represents a nontrivial class of $H^1_{{\rm loc}} ( \Gal ( k ( \A[p^n] ) / k ) , \A[p^n] )$. 
Conclude by applying Proposition \ref{profin} with $Q$ in the place of $D$ and $W$ in the place of $Z$.
\CVD   

Before finding the counterexample for $p = 2$ to Theorem \ref{teo2}, we need the following result on the Galois action over the $p$-torsion of an elliptic curve.

\begin{lem}\label{lem32}
Given a prime number $p$, a positive integer $n$ and a subgroup $G$ of ${\rm GL}_2 ( \Z / p^n \Z )$, there exists a number field $k$ and an elliptic curve $\E$ defined over $k$ such that there are an isomorphism $\phi \colon {\rm Gal} ( k ( \E[p^n] ) / k ) \rightarrow G$ and a $\Z / p^n \Z$-linear isomorphism $\tau \colon \E[p^n] \rightarrow ( \Z / p^n \Z )^2$ such that 
\begin{equation}\label{eqn:rel3ancora}
\phi ( \sigma ) \cdot \tau ( v ) = \tau ( \sigma ( v ) ).
\end{equation}
for all $\sigma \in \Gal ( k ( \E[p^n] ) / k )$ and $v \in \E[p^n]$.
\end{lem}    

\DIM Greicius \cite{G} found an example of an elliptic curve defined over a cubic field such that the morphism of the absolute Galois group into ${\rm GL}_2 ( \hat{Z} )$ induced by the Galois action over the torsion points, is surjective. 
Zywina showed that this is the case for most elliptic curves defined over number fields not intersecting cyclotomic fields (see \cite{Z}).
Then there exists a number field $L$ and an elliptic curve $\E$ defined over $L$ such that $\Gal ( L ( \E[p^n] ) / L )$ is isomorphic to the all group of the $\Z / p^n \Z$-linear automorphisms of $\E[p^n]$. 
Such isomorphism may be described in the following way : 
let $\tau$ be a $\Z / p^n \Z$-linear isomomorphism of $\E[p^n]$ to $( \Z / p^n \Z )^2$ and let $P_1$ respectively $P_2$ be the elements of $\E[p^n]$ such that $\tau ( P_1 ) = ( 1, 0 )$, respectively $\tau ( P_2 ) = ( 0, 1 )$. 
For every $\sigma$ in ${\rm Gal} ( L ( \E[p^n] ) / L )$, there exist $a_\sigma$, $b_\sigma$, $c_\sigma$, $d_\sigma$ in $\Z / p^n \Z$ such that $\sigma ( P_1 ) = a_\sigma P_1 + c_\sigma P_2$ and $\sigma ( P_2 ) = b_\sigma P_1 + d_\sigma P_2$.
Define $\phi \colon {\rm Gal} ( L ( \E[p^n] ) / L ) \rightarrow {\rm GL}_2 ( \Z / p^n \Z )$ the function sending $\sigma \in {\rm Gal} ( L ( \E[p^n] ) / L )$ to the matrix 
\[
\phi ( \sigma ) = \left(
\begin{array}{cc}
a_\sigma & b_\sigma \\
c_\sigma & d_\sigma  \\
\end{array}
\right).
\]
An easy verification shows that $\phi$ is an injective homomorphism and that for every $\sigma \in {\rm Gal} ( L ( \E[p^n] ) / L )$ and $v \in \E[p^n]$, we have $\phi ( \sigma ) \cdot \tau ( v ) = \tau ( \sigma ( v ) )$.
Moreover, since $\Gal ( L ( \E[p^n] ) / L )$ is isomorphic to the all group of the $\Z / p^n \Z$-linear automorphisms of $\E[p^n]$, $\phi$ is an isomorphism.

Let $k$ be the subfield of $L ( \A[p^n] )$ fixed by $\phi^{-1} ( G )$.
Then by Galois correspondence, $\Gal ( k ( \E[p^n] ) / k ) = \phi^{-1} ( G )$.
Thus, the restriction of $\phi$ to $\Gal ( k ( \E[p^n] ) / k )$ gives an isomorphism between $\Gal ( k ( \E[p^n] ) / k )$ and $G$ that satisfies the relation (\ref{eqn:rel3ancora}).
\CVD  

\begin{rem}\label{rem33}
To prove that Theorem \ref{teo2} is false for $p=2$ we just need to prove Lemma \ref{lem32} in the particular case $p=2$.
Then, it is sufficient to prove that there exists a number field $L$ and an elliptic curve defined over $L$ with surjective $2$-adic representation. 
We thank the referee to point out the very interesting paper \cite{JR}.
There, see \cite[example 5.4]{JR}, the authors prove that the elliptic curve over $\Q$ given by the equation $y^2 + y = x^3 - x$ has surjective $2$-adic representation.
\end{rem}   

\begin{pro}\label{pro34}
Theorem \ref{teo2} is false for $p = 2$.
\end{pro}

\DIM Let $G$ be the group $( \Z / 8 \Z )^\ast$ acting in the classical way over $\Z / 8 \Z$.
Then $H^1_{{\rm loc}} ( G , \Z / 8 \Z ) \neq 0$ (see also \cite{DZ}).
In fact let $b$ be a generator of $\Z / 8 \Z$, and for $i \in \{1 , 3 , 5 , 7\}$, call $g_i \in G$ such that $g_i b = i b$. 
The cocycle $Z \colon G \rightarrow \Z / 8 \Z$ sending $g_1$ and $g_7$ to $0$, and $g_3$ and $g_5$ to $4 b$ is not a coboundary and its class is in $H^1_{{\rm loc}} ( G , \Z / 8 \Z )$. 

In order to apply corollary \ref{cor31} we will work in $(\Z/16\Z)^2$ as the point $Q$ that we will define will be a point of order $16$ of some elliptic curve $\E$ defined over a number field $k$. 
Let $G_{16} \leq {\rm GL}_2 ( \Z / 16 \Z )$ be the group 
\[
\bigg \langle
\left(
\begin{array}{cc}
3 & 8 \\
0 & 1  \\
\end{array}
\right),
\left(
\begin{array}{cc}
15 & 0 \\
0 & 1  \\
\end{array}
\right)
\bigg \rangle,
\]
and $G_8 \leq {\rm GL}_2 ( \Z / 8 \Z )$ be the image of $G_{16}$ by the homomorphism induced by the projection $\Z / 16 \Z \rightarrow \Z / 8 \Z$.
Then $G_8$ is the group of
\[
\bigg\{
\left(
\begin{array}{cc}
a & 0 \\
0 & 1  \\
\end{array}
\right),  \ a \in ( \Z / 8 \Z )^\ast
\bigg\}.
\]
By Lemma \ref{lem32} (see also Remark \ref{rem33}), there exist a number field $k$ and an elliptic curve $\E$ defined over $k$ such that there exist an isomorphism $\phi_{16} \colon \Gal ( k ( \E[16] ) / k ) \rightarrow G_{16}$ and a $\Z / 16 \Z$-linear isomorphism $\tau_{16} \colon \E[16] \rightarrow ( \Z / 16 \Z )^2$ such that $\phi_{16} ( \sigma ) \cdot \tau_{16} ( v ) = \tau_{16} ( \sigma ( v ) )$, for every $\sigma \in \Gal ( k ( \E[16] ) / k )$ and $v \in \E[16]$. 
Let $\pi_1 \colon \Gal ( k ( \E[16] ) / k ) \rightarrow  \Gal ( k ( \E[8] ) / k )$ be the morphism sending $\sigma \in \Gal ( k ( \E[16] ) / k )$ to its restriction to $k ( \E[8] )$, let $\pi_2 \colon G_{16} \rightarrow G_8$ be the homomorphism induced by the projection $\Z / 16 \Z \rightarrow \Z / 8 \Z$, let $\pi_3 \colon \E[16] \rightarrow \E[8]$ be the multiplication by $2$ and let $\pi_4 \colon ( \Z / 16 \Z )^2 \rightarrow ( \Z / 8 \Z )^2$ be the projection.
It is straightforward to prove that there exist unique $\phi_8 \colon \Gal ( k ( \E[8] ) / k ) \rightarrow G_8$ and $\tau_8 \colon \E[8] \rightarrow ( \Z / 8 \Z )^2$ making the following diagrams commutative
\[
\xymatrix{
\Gal ( k ( \E[16] ) / k ) \ar[r]^-{\phi_{16}} \ar[d]_{\pi_1} & G_{16} \ar[d]^{\pi_2}\\
\Gal ( k ( \E[8] ) / k ) \ar[r]^-{\phi_8} & G_8
} 
\]

\[
\xymatrix{
\E[16] \ar[r]^-{\tau_{16}} \ar[d]_{\pi_3}  & ( \Z / 16 \Z )^2 \ar[d]^{\pi_4} \\
\E[8] \ar[r]^-{\tau_8} & ( \Z / 8 \Z )^2
}
\]

Moreover, for all $\gamma \in \Gal ( k ( \E[8] ) / k )$ and $w \in \E[8]$, we have \\ $\phi_8 ( \gamma ) \cdot  \tau_8 ( w ) = \tau_8 ( \gamma ( w ) )$.

Let $Q_1$ be $\tau_{16}^{-1} ( ( 1, 0 ) )$, let $Q_2$ be $\tau_{16}^{-1} ( ( 0, 1 ) )$.
Moreover set $P_1 = 2 Q_1$ and $P_2 = 2 Q_2$.
Finally, let $\widetilde{\tau_8}$ be the homomorphism from $\langle P_1 \rangle$ to $\Z / 8 \Z$ sendind $P_1$ to $1$ and let $\widetilde{\phi_8}$ be the homomorphism from $\Gal ( k ( \E[8] ) / k )$ to $( \Z / 8 \Z )^\ast$ sending $\Gal ( k ( \E[8] ) / k )$ to $a_\sigma$, where $\sigma ( P_1 ) = a_\sigma P_1$. 
Then, if we set $p = 2$, $n = 3$, $d=2$, $r = 1$, $\A = \E$, $\phi = \widetilde{\phi_8}$ and $\tau = \widetilde{\tau_8}$, all the hypotheses of Corollary \ref{cor31} until the point 2. hold.
Then to apply Corollary \ref{cor31} and to conclude to proof, it is sufficient to find a point $Q \in \E ( k ( \E[8] ) )$ such that $Q$ has order a power of $2$ and, for every $\sigma \in \Gal ( k ( \E[8] ) / k )$, we have $\sigma ( Q ) = Q + \tau_8^{-1} ( Z_{\phi ( \sigma )} )$.
Let us prove that $Q_2$ is such a point.
By definition $Q_2$ has order $16$ and so its order is a $2$ power.
Clearly $Q \in \E ( k ( \E[16] ) )$.
Observe that $\Gal ( k ( \E[16] ) / k ( \E[8] ) )$ is the inverse image by the isomorphism $\phi_{16}$ of the subgroup $H$ of $G_{16}$ fixing the group $( 2 \Z / 16 \Z )^2$.
Then $H$ is the group generated by 
\[
\bigg \langle
\left(
\begin{array}{cc}
15 & 0 \\
0 & 1 
\end{array}
\right)
\bigg \rangle.
\]
Thus $H$ fixes $( 0, 1 )$ and so $\phi_{16}^{-1} ( H ) = \Gal ( k ( \E[16] ) / k ( \E[8] ) )$ fixes $\tau_
{16}^{-1} ( ( 0, 1 ) ) = Q_2$.
Thus $Q_2 \in \E ( k ( \E[8] ) )$.
Finally we should prove that for all $\sigma \in \in \Gal ( k ( \E[8] ) / k )$, we have $\sigma ( Q_2 ) = Q_2 + \widetilde{\tau_8}^{-1} ( Z_{\widetilde{\phi_8} ( \sigma )} )$, but this is an easy verification.
\CVD    

\section{The case of a ${\rm GL}_2$-type variety}\label{sec4}

Recall that an abelian variety $A$ defined over a number field $k$ is said to be of ${\rm GL}_2$-type if there exists a number field $E$ such that $[E: \Q] = {\rm dim} ( \A )$ and an embedding $\phi \colon E \hookrightarrow {\rm End}_k ( \A ) \otimes \Q$ (see \cite[Section 2, Section 5]{R}). We note $\O_E$ the endomorphism ring of $E$. Then $R = E \cap \phi^{-1} \left( {\rm End}_k ( \A ) \otimes \Z \right)$ is an order of $\O_E$ and $\phi$ induce an embedding $R \hookrightarrow {\rm End}_k ( \A )$. 

Following \cite{R}, we say that a prime number is {\bf good} for $\A$ if it does not divide the index $[\O_E : R]$. (In fact the notion of good depends also of the choice of $E$ if there exists more than a number field holding this property). 

Ribet \cite[Proposition 2.2.1]{R} proved the following result.

\begin{pro}\label{pro41}
If $p$ is good for $\A$, then the Tate module $T_p ( \A )$ is a free $R_p = R \otimes \Z_p$-module of rank $2$; equivalently, it is a free $\O_p = \O_E \otimes \Z_p$-module of rank $2$.
\end{pro}

Then we have the following corollary.

\begin{cor}\label{cor42}
Let $p$ be a good prime for $\A$. Then for every positive integer $n$, $\Gal ( k ( \A[p^n] ) / k )$ is isomorphic to a subgroup of $\prod_{\P \mid p} \GL_2 ( \O_E / \P^{e_\P n} )$, where $\P$ is a prime of $R$ dividing $p$ and $e_\P$ is the ramification index of $\P$.
\end{cor}

\DIM For every $\P$ dividing a good prime $p$ and every $m \in \N$, let $\A[\P^m]$ be the set of the elements of $\A  ( \overline{k} )$ killed by $\P^m$.  
By Proposition \ref{pro41}, $\A[p^m]$ is isomorphic to $( \O_E / p^m \O_E )^2$.
Observe (or see \cite{R}) that $\A[\P^m] \simeq ( \O_E / \P^{e_\P m} \O_E )^2$ and $\A[p^m] \simeq \oplus_{\P \mid p} \A[\P^{e_\P m}]$. 
Since we have an injective group homomorphism 
\[
\phi_p \colon \Gal ( k ( \A[p^n] ) / k ) \rightarrow {\rm Aut} ( \A[p^n] ) \simeq \prod_{\P \mid p} {\rm Aut} ( \A[\P^{e_\P n}] )
\]
and ${\rm Aut} ( \A[\P^{e_\P n}] )$ is isomorphic to a subgroup of $\GL_2 ( \O_E / \P^{e_\P n} )$, we conclude by applying the Chinese remainder Theorem.
\CVD

Let $p$ be a good prime number. 
Denote $\phi_{\P}$ the composition between $\phi_p$ (see above the proof of Corollary \ref{cor42}) and the projection
$\prod_{\P^\prime \mid p} {\rm Aut} ( \A[\P^{\prime{e_\P n}}] ) \rightarrow {\rm Aut} ( \A[\P^{e_\P n}] )$. 
From now on, let us simplify our notation by putting $G_{p, n} = \Gal ( k ( \A[p^n] ) / k )$ and $G_{\P, n} = \phi_{\P} ( G_{p, n} )$.
 
The following lemma tells us that, studying the local cohomology of $\A[p^n]$ is equivalent to studying the local cohomology of the $G_{\P , n}$-modules $\A[\P^{e_\P n}]$.

\begin{lem}\label{lem43}
Let $p$ be a good prime for $\A$ and, for every prime ideal $\P$ in $R$, let $e_\P$ be its ramification index.
Then there is an isomorphism
\[
H^1_{{\rm loc}} ( G_{p , n} , \A[p^n] ) \rightarrow \bigoplus_{\P \mid p}  H^1_{{\rm loc}} ( G_{\P , n}, \A[\P^{e_\P n}] ). 
\]   
\end{lem}

\DIM Recall that
\[
\A[p^n] \simeq \bigoplus_{\P \mid p} ( \A[\P^{e_\P n}] )
\] 
as $G_{p, n}$-module.
Then there is an isomorphism (see for instance \cite{S})
\[
H^1 ( G_{p, n} , \A[p^n] ) \rightarrow \bigoplus_{\P \mid p}  H^1 ( G_{p, n} , ( \A[\P^{e_\P n}] )
\]
that induces an isomorphism
\[
H^1_{{\rm loc}} ( G_{p, n} , \A[p^n] ) \rightarrow \bigoplus_{\P \mid p}  H^1_{{\rm loc}} ( G_{p, n} , \A[\P^{e_\P n}] ).
\]
Then it suffices to show that for every $\P$ over $( p )$ we have an isomorphism 
\[
H^1_{{\rm loc}} ( G_{\P , n} , ( \A[\P^{e_\P n}] ) \simeq H^1_{{\rm loc}} ( ( G_{p, n} , \A[\P^{e_\P n}] ).
\]
The inflation-restriction sequence (see \cite[Section VII, Proposition 4]{S}) that we obtain by choosing $\ker ( \phi_{\P} )$ as normal subgroup of $G_{p, n}$ is
{\small
\[
0 \rightarrow H^1 ( G_{p, n} / \ker ( \phi_{\P} ), ( \A[\P^{e_\P n}] ) \rightarrow H^1 ( G_{p, n} , \A[\P^{e_\P n}] ) \rightarrow H^1 ( \ker ( \phi_{\P} ) , \A[\P^{e_\P n}] ).  
\]
}
Since $G_{p, n} / \ker ( \phi_{\P} )$ is isomorphic to $G_{\P , n}$, the inflation induces an injective function 
\[
H^1 ( G_{\P , n}, \A[\P^{e_\P n}] ) \rightarrow H^1 ( G_{p , n} , \A[\P^{e_\P n}] ).
\]
Moreover $\ker ( \phi_{\P} )$ acts like the identity over $\A[\P^{e_\P n}]$.
Then, the same argument in \cite[p. 322]{DZ} for proving their Proposition 2.1 shows that the inflation induces an isomorphism 
\[
H^1_{{\rm loc}} ( G_{\P , n}, ( \A[\P^{e_\P n}] ) \simeq H^1_{{\rm loc}} ( G_{p , n} , \A[\P^{e_\P n}] ).
\]
This concludes the proof.
\CVD  

\begin{rem}\label{rem44}
We make more explicit the isomorphism 
\[
H^1_{{\rm loc}} ( G_{p , n} , \A[p^n] ) \rightarrow \bigoplus_{\P \mid p}  H^1_{{\rm loc}} ( G_{\P , n}, \A[\P^{e_\P n}] )
\]
in Lemma \ref{lem43}, because it will be useful in the proof of Theorem \ref{teo3}.
Let $Z$ be a cocycle representing a class $[Z]$ of $H^1_{{\rm loc}} ( G_{p, n} , \A[p^n]$.
Let $\pi_\P$ be the projection of $\A[p^n]$ to $\A[\P^{e_\P n}]$.
Then the image of $[Z]$ is the sum of the classes representing the cocycles $Z^\P \colon G_{\P, n} \rightarrow \A[\P^n]$ such that, for every $\sigma \in G_{p , n}$, $Z^\P_{\phi_\P ( \sigma_\P )} = \pi_\P ( Z_\sigma )$.
\end{rem}

\subsection{Proof of Theorem \ref{teo3} and a particular case in dimension 2.}

In this section we prove Theorem \ref{teo3}. After that we show that in the special case of $\GL_2$-type abelian surfaces $\A$ defined over a number field $k$ such that ${\rm End}_{\overline{k}} ( \A )$ is an order in a quadratic real field $E$, and that all endomorphisms of $\A$ are defined over $k$, the hypotesis on the inertia degree in this theorem is necessary, for $p$ big enough compared to $\A$. Observe that such $\A$ is a ${\rm GL}_2$-type variety according to our definition. 
In fact, if all endomorphisms of $\A$ are defined over $k$, we have ${\rm End}_{\overline{k}} ( \A ) = {\rm End}_k ( \A )$.
Since ${\rm End}_k ( \A )$ is an order in the totally real quadratic field $E$, we have that $E$ embeds into ${\rm End}_k ( \A ) \otimes \Q$, and so ${\rm End}_k ( \A ) \otimes \Q$ contains a number field of degree equal to the dimension of $\A$. 
Thus $\A$ is a ${\rm GL}_2$-type variety.

 \vspace{5 pt}

{\bf Proof of Theorem \ref{teo3}.} 
Let $p$ be an odd good prime for $\A$ that splits completely over $\O_E$.
Suppose that there exists a torsion point that is a counterexample for the local global divisibility by a power of $p$ over $\A ( k )$. Then by Lemma \ref{lem11} there exists $P \in \A ( k )$ of order a power of $p$ such that the local-global divisibility by a power of $p$ fails for $P$ over $\A ( k )$. 
Let $p^n$ be the minimal power of $p$ for which this occurs. Then, by Lemma \ref{lem12}, there exists $Q \in \A ( k )$ such that $p^{n-1} Q = P$. In particular, $\A$ admits a point of order $p^n$ defined over $k$ of the shape $p^aQ$.
Let $D \in \A ( \overline{k} )$ be such that $p D = Q$ (and so $p^n D = P$).
Then, by \cite[Corollary 2.3]{DZ}, $D \in \A ( k ( \A[p^n] ) )$ and the cocycle $Z \colon G_{p, n} \rightarrow \A[p^n]$ sending $\sigma \in G_{p, n}$ to $Z_\sigma = \sigma ( D ) - D$ represents a non-trivial class $[Z]$ of $H^1_{{\rm loc}} ( G_{p, n} , \A[p^n] )$ (see the proof of \cite[Proposition 2.1]{DZ}).
Write $P = \sum_{\P \mid p} P_\P$, $Q = \sum_{\P \mid p} Q_\P$ and $D = \sum_{\P \mid p} D_\P$ with $P_\P \in \A[\P^b]$, $Q_\P \in \A[\P^{b+n-1}]$ and $D_\P \in \A[\P^{b+n}]$ for a certain positive integer $b$.
By Lemma \ref{lem43}, we have an isomorphism
\[
H^1_{{\rm loc}} ( G_{p , n} , \A[p^n] ) \rightarrow \bigoplus_{\P \mid p}  H^1_{{\rm loc}} ( G_{\P , n}, \A[\P^n] ).
\]
Denote as before $\phi_\P$ the projection of $G_{p, n}$ to $G_{\P , n}$ and $\pi_\P$ the projection of $\A[p^n]$ to $\A[\P^n]$. 
Recall, see Remark \ref{rem44}, that the image of $[Z]$ for this ismorphism is the sum of the classes representing the cocycles $Z^\P \colon G_{\P, n} \rightarrow \A[\P^n]$ such that, for every $\sigma \in G_{p , n}$, $Z^\P_{\phi_\P ( \sigma_\P )} = \pi_\P ( Z_\sigma )$.
In particular, since $Z_\sigma = \sigma ( D ) - D$ for every $\sigma \in G_{p, n}$, we have $Z^\P_{\phi_\P ( \sigma )} = \phi_\P ( \sigma ) ( D_\P ) - D_\P$.  

Let $\P$ be such that $P_\P = 0$.
Then $D_\P \in \A[\P^n]$ and so $Z^P$ is the coboundary associated to $D_\P$.
Thus $[Z]$ is in the kernel of the projection $H^1_{{\rm loc}} ( G_{p , n} , \A[p^n] ) \rightarrow \bigoplus_{P_\P = 0}  H^1_{{\rm loc}} ( G_{\P , n}, \A[\P^n] )$.

Let $\P$ be such that $P_\P \neq 0$.
Since $p D = Q \in \A ( k )$ and $p^{n-1} Q = P$, in particular $p^{n-1} Q_\P = P_\P$ and so $G_{\P , n}$ fixes an element of $\A[\P^n]$ of order $\geq p^n$.
Since $( p )$ splits completely in $\O_E$, we have that $A[\P^n]$ is isomorphic to $( \Z / p^n \Z )^2$ via an ismorphism $\tau$, $G_{\P , n}$ is isomorphic to a subgroup of ${\rm GL}_2 ( \Z / p^n \Z )$ fixing an element of order $p^n$ of $( \Z / p^n \Z )^2$ via an isomorphism $\phi$, and for every $v \in \A[\P^n]$ and $\sigma \in G_{\P, n}$, $\phi ( \sigma ) \cdot \tau ( v ) = \tau ( \sigma ( v ) )$.
Then by Proposition \ref{pro21}, we get $H^1_{{\rm loc}} ( G_{\P , n}, \A[\P^n] ) = 0$.
Thus $Z^\P$ is a coboundary for every $\P$.
But this contradicts the fact that $[Z] \neq 0$.
\CVD

We now show that, in the particular case described above and $p$ big enough compared to $\A$, the hypothesis that $( p )$ has inertia degree $1$ is necessary in Theorem \ref{teo3}.
We first need the following lemma.

\begin{lem}\label{lem51}
Let $p$ be a prime number, let $q$ be $p^2$ and let $\alpha$ be a generator of $\F_q$.  
Let $G$ be the subgroup of ${\rm GL}_2 ( \F_q )$, such that
\[
G = \bigg \{ \sigma ( \lambda_1 , \lambda_2 ) =
\left(
\begin{array}{cc}
1 & \lambda_1 + \lambda_2 \alpha \\
0 & 1 \\
\end{array}
\right), 
\ \lambda_1 , \ \lambda_2 \in \F_p 
\bigg \},
\]
acting over $\F_q^2$ in the classical way. Then the function $Z \colon G \rightarrow \F_q^2$ sending $\sigma ( \lambda_1 , \lambda_2 )$ to $( \lambda_2 , 0 )$ is a cocycle and it represents a non-trivial class of $H^1_{{\rm loc}} ( G , \F_q^2 )$.
\end{lem}

\DIM First observe that $Z$ is a cocycle because it is a group homomorphism and $G_1$ fixes the image of $Z$.

Let $( \lambda_1 , \lambda_2 )$ be in $( \F_p )^2$ with $\lambda_2 \neq 0$.
There exists $\beta_{( \lambda_1 , \lambda_2 )} \in \F_q$ such that $( \lambda_1 + \lambda_2 \alpha ) \beta_{( \lambda_1 , \lambda_2 )} = \lambda_2$.
Then $Z_{\sigma ( \lambda_1 , \lambda_2 )} = \sigma ( \lambda_1 , \lambda_2 ) (  0 , \beta_{( \lambda_1 , \lambda_2 )} ) - ( 0, \beta_{( \lambda_1 , \lambda_2 )} )$.
On the other hand if $\lambda_2 = 0$, we have $Z_{\sigma ( \lambda_1 , \lambda_2 )} = \sigma ( \lambda_1 , \lambda_2 ) (  0 , 0 ) - ( 0, 0 )$.
This proves that $Z$ satisfies the local conditions.

Finally observe that for every $( x, y ) \in \F_q^2$ such that $Z_{\sigma ( 0 , 1 )} = \sigma ( 0 , 1 ) ( x , y  ) - ( x , y )$, we have $y \neq 0$ and, for every $( w, t ) \in ( \F_{p^2} )^2$ such that $Z_{\sigma ( 1 , 0 )} = \sigma ( 1 , 0 ) ( w , t  ) - ( w , t )$, we have $t = 0$.
Then $Z$ is not a coboundary.
This proves that $Z$ represents a non-trivial class in $H^1_{{\rm loc}} ( G , \F_q^2 )$.
\CVD 

We also need the following result on the Galois action over the $p$-torsion points of a ${\rm GL}_2$-type surface over $\Q$.

\begin{lem}\label{lem52}
Let $\A$ be an abelian surface defined over a number field $k$ such that ${\rm End}_{\overline{k}} ( \A )$ is an order in a quadratic real field $E$, and that all endomorphisms of $\A$ are defined over $k$, then there exists 
a number field $L$ containing $k$ and a constant $C$, such that for every prime number $p > C$, for every positive integer $n$ there is an isomorphism $\phi_n \colon {\rm Gal} ( L ( \A[p^n] ) / L ) \rightarrow G_n$, where $G_n$ is the subgroup of ${\rm GL}_2 ( \O_E / p^n \O_E )$ of the elements with determinant in $( \Z / p^n \Z )^\ast$, and a $\O_E / p^n \O_E$-linear isomorphism $\tau_n \colon \A[p^n] \rightarrow ( \O_E / p^n \O_E )^2$ such that 
\begin{equation}\label{eqn:rel53}
\phi_n ( \sigma ) \cdot \tau_n ( v ) = \tau_n ( \sigma ( v ) )
\end{equation}
for all $\sigma \in \Gal ( L ( \A[p^n] ) / L )$ and $v \in \A[p^n]$.
\end{lem}

\DIM Lombardo \cite[Theorem 1.4, Corollary 1.5]{L} proved that if $\A$ is an abelian surface defined over a number field $k$ such that ${\rm End}_{\overline{k}} ( \A )$ is an order in a quadratic real field $E$, and that all endomorphisms of $\A$ are defined over $k$, then there exists a constant $C$ such that for every prime number $p > C$, for every positive integer $n$ there is an isomorphism $\phi_n \colon {\rm Gal} ( L ( \A[p^n] ) / L ) \rightarrow G_n$, where $G_n$ is the subgroup of ${\rm GL}_2 ( \O_E / p^n \O_E )$ of the elements with determinant in $( \Z / p^n \Z )^\ast$.
The proof of the existence of the morphism $\tau_n$ and the relation (\ref{eqn:rel53}) is now very similar  to the proof of Lemma \ref{lem32}, and we omit it.
\CVD

We can now show that for this class of $\GL_2$-type abelian varieties, the hypothesis on the inertia degree in our Theorem \ref{teo3} is necesary. 

\begin{pro}\label{pro54}
Let $\A$ be an abelian surface defined over a number field $k$ such that ${\rm End}_{\overline{k}} ( \A )$ is an order in a  quadratic real field $E$, and that all endomorphisms of $\A$ are defined over $k$,
There exists a good prime number $p$ inert in $\O_E$, such that there exists a point $P \in \A ( L )$ of order $p$, which is divisible by $p$ for all but finitely many valuations $v$ of $L$, but it is not divisible by $p$ over $\A ( L )$.
\end{pro}

\DIM Recall that given a number field $F$ the set of the prime numbers inert in $F$ is infinite and that the set of bad primes of a ${\rm GL}_2$-type variety is finite.
Then by Lemma \ref{lem52}, we can choose a ${\rm GL}_2$-type surface $\A$ defined over a number field $k$ such that a real quadratic number field $E$ embeds into ${\rm End}_k ( \A ) \otimes \Q$, and a prime number $p$ inert over $E$ such that there are isomorphisms $\phi_n \colon {\rm Gal} ( k ( \A[p^n] ) / k ) \rightarrow G_n$, where $G_n$ is the subgroup of ${\rm GL}_2 ( \O_E / p^n \O_E )$ of the elements with determinant in $( \Z / p^n \Z )^\ast$, and a $\O_E / p^n \O_E$-linear isomorphism $\tau_n \colon \A[p^n] \rightarrow ( \O_E / p^n \O_E )^2$ such that 
\[
\phi_n ( \sigma ) \cdot \tau_n ( v ) = \tau_n ( \sigma ( v ) )
\]
for all $\sigma \in \Gal ( \Q ( \A[p^n] ) / \Q )$ and $v \in \A[p^n]$.

Since $( p )$ is inert, $\O_E / p \O_E$ is isomorphic to $\F_q$, where $q = p^2$.  
Let $\alpha$ be a generator of $\F_q$. 
Consider the subgroup $H_1$ of $G_1$
\[
H_1 = \bigg \{ \sigma ( \lambda_1 , \lambda_2 ) =
\left(
\begin{array}{cc}
1 & \lambda_1 + \lambda_2 \alpha \\
0 & 1 \\
\end{array}
\right), 
\ \lambda_1 , \ \lambda_2 \in \F_p 
\bigg \}.
\]
By Lemma \ref{lem51} the function $Z \colon H_1 \rightarrow \F_q^2$ sending $\sigma ( \lambda_1 , \lambda_2 )$ to $( \lambda_2 , 0 )$ is a cocycle and it represents a non-trivial class of $H^1_{{\rm loc}} ( G , \F_q^2 )$.
Consider the injective group homomorphism $\psi \colon ( \O_E / p \O_E ) \rightarrow ( \O_E / p^2 \O_E ) )$, sending $x \in \O_E$ modulo $( p )$ in the class $p x$ modulo $( p^2 )$.     
Then define    
\[      
H_2 = \bigg \{ \sigma^\prime ( \mu_1 , \mu_2 ) =
\left(
\begin{array}{cc}
1 +  p \mu_2 & \mu_1 + \mu_2 \alpha^\prime \\
0 & 1 \\
\end{array}
\right), 
\ \mu_1 , \ \mu_2 \in \Z / p^2 \Z  
\bigg \},
\]
where $\alpha'\in O_E/p^2O_E$ is such that $\psi ( \alpha ) = p \alpha^\prime$.
Then $H_2$ is a subgroup of $G_2$.
Using the notation of Lemma \ref{lem52}, consider the subfield $L$ of $k ( \A[p^2] )$ fixed by the group $\phi^{-1} ( H_2 )$.
Moreover let $D_1$ be $\tau_2^{-1} ( ( 1, 0 ) )$ and $D_2 = \tau_2^{-1} ( ( 0, 1 ) )$.
We shall prove that $p D_1$ is in $\A ( L )$, is divisible by $p$ over $\A ( L_v )$ for all but finitely many valuations $v$ of $L$, but is not divisible by $p$ over $\A ( L )$.
First of all observe that $\Gal ( L ( \A[p] ) / L ) = \phi_1^{-1} ( H_1 )$, and so $p D_1 = \tau_2^{-1} ( ( p, 0 ) ) = \tau_1^{-1} ( ( 1, 0 ) )$ is clearly fixed by $\tau^{-1} ( H_1 )$.
Thus it is in $\A ( L )$.
Let us observe that $D_1 \in \A ( L ( \A[p] ) )$. 
In fact let $\gamma \in \phi_2^{-1} ( H_2 )$ fixing $\A[p]$. 
Thus $\phi_2 ( \gamma ) = \sigma^\prime ( \mu_1 , \mu_2 )$ for some $\mu_1 \in p \Z / p^2 \Z$ and $\mu_2 \in p \Z / p^2 \Z$.
Thus $\sigma^\prime ( ( \mu_1 , \mu_2 ) ) ( 1, 0 ) = ( 1, 0 )$ and so $D_1 = \tau_2^{-1} ( ( 1, 0 ) )$ is fixed by $\gamma$. 
Let $W \colon \Gal ( L ( \A[p] ) / L ) \rightarrow \A[p]$ the cocycle defined by $\tau_1^{-1} \circ Z \circ \phi_1$. 
It is easy to verify that since $Z$ represents a non-trivial class of $H^1_{{\rm loc}} ( H_1 , \F_q^2 )$, then $W$ represents a non-trivial class of $H^1_{{\rm loc}} ( \Gal ( L ( \A[p] ) / L ) , \A[p] )$.
Moreover another easy verification shows that for every $\delta \in \Gal ( L ( \A[p] ) / L )$, we have $\delta ( D_1 ) = D_1 + Z_\delta$.
Thus by Proposition \ref{profin}, $p D_1$ gives a counterexample to the local-global divisibility principle by $p$ over $\A ( L )$.
\CVD

\begin{rem}\label{rem55}
We try to partially answer to the following two questions: let $\A$ be a ${\rm GL}_2$-type variety defined over a number field $k$ such that the number field $E$ of degree equal to the dimension of $\A$ embeds into ${\rm End}_k ( \A ) \otimes \Q$.
Does the local-global divisibility principle by a power of a good ramified prime over $E$ holds for $\A ( k )$ ? Does the local-global divisibility principle by a power of a bad prime holds for $\A ( k )$?

Suppose that $( p )$ is a good prime ramified over $\O_E$.
In this case Corollay \ref{cor42} applies and so the image of the Galois representation is isomorphic to a subgroup of ${\rm GL}_2 ( \O_E / ( p^n ) )$.
On the other hand to our knowledge there does not exist any results on the Galois image over $\A[p]$ in the case when $p$ is ramified (then we do not have a lemma analogous to Lemma \ref{lem52}).

Suppose now that $p$ is a bad prime. 
In this case Corollary \ref{cor42} does not apply and we think we should study the representation of the absolute Galois group over ${\rm GL}_{2d} ( \Z / p^n \Z )$, where $d$ is the dimension of $\A$.
Then the problem is as difficult as studying the local-global divisibility problem for an arbitrary abelian variety of dimension $d$ and to our knowledge we do not have any result analogous to Lemma \ref{lem52} in this case.
\end{rem}

\section{Proof of Corollary \ref{cor4}}\label{sec5}  
 
Let $p\neq 2$ be a prime number wich satisfies the conditions of the corollary. We shall show that there exists an order $R$ in a quadratic field $E$ such that $R$ is optimally embedded in $\O$ and $p$ is a good prime for $R$ that split in $E$. 
Observe that the condition that $B$ is a quaternions algebra over $\Q$ indefinite implies that the maximal orders of $B$ are all in the same conjugacy class. For $d\in \N$ such that $d$ is not a square, let $E=\Q(\sqrt{d})$ be a quadratic field. Let $\O_E$ be the integer rings of $E$. By a theorem of Eichler there exists a maximal order $\O_m'$ of $B_D$ such that $\O_E$ embeds maximally in $\O_m'$ if and only if all prime $\ell$ dividing the discriminant of $B$, $\ell$ does not split in $E$  (cf. for instance \cite[thm 5.11 and cor 5.12, pp. 92 and 94]{V} for the number of such embeddings). An optimal embedding from $\O_E$ to $\O_m'$ is an embedding $\varphi\colon E\rightarrow B$ such that $\varphi^{-1}(\O_m')=O_E$ (or equivalently $\varphi(R)=O_m'\cap\varphi(E)$). As we said all maximal orders are in the same conjugacy class, so if $\O_E$ is optimally embedded in $\O_m'$ it is also optimally embedded in $\O_m$. By the Chinese remainder theorem there exists a solution not square $d$ to the system : 
$$
\begin{cases}
\left(\topbot{d}{p}\right)=1  \\
\left(\topbot{d}{\ell}\right)=-1 \ \text{ for all } \ell\mid D \ \ \ell\neq 2  \\
d = 5 \mod 8 \ \ \text{ if } 2\mid D
\end{cases}
$$
For this $d$ by Eichler theorem, we have an embedding $\varphi\colon E=\Q(\sqrt{ d }) \rightarrow B$ that is an optimal embedding from $\O_E$ to $\O_m$. Recall that by hypothesis we have : $[\O_m:\O]=M$ so $M\O_m\subseteq \O\subseteq \O_m$. Let $R=\varphi^{ -1 } (\O)$ then $R$ is an order of $E$ and we have $M\O_E=M\varphi^{ -1 } (\O_m) \subseteq \varphi^{ -1 } (\O)=R $. So the condition $p\nmid M$ implies that $p$ does not divide the conductor of $R$. 

We have shown that there exists a quadratic number field $E$ and an embedding $\varphi\colon E\rightarrow B \simeq {\rm End}_k(\A) \otimes \Q$, such that $p$ split in $E$ and $\varphi^{-1}({\rm End}_k(\A)) = R$ is an order of $E$ such that $p$ does not divide $[\O_E:R]$. We may apply Theorem \ref{teo3}, so for all $n\geq 1$ the local global divisibility by $p^n$ holds for the torsion points of $\A$. 
\CVD


\end{document}